\documentclass{amsart}

\usepackage{fancyhdr}
\usepackage{lastpage}
\usepackage{stmaryrd,yhmath}

\newcommand{\bdis}{\begin{displaymath}}
\newcommand{\edis}{\end{displaymath}}
\newcommand{\be}{\begin{equation}}
\newcommand{\ee}{\end{equation}}
\newcommand{\mbb}{\mathbb}
\newcommand{\mcal}{\mathcal}

\newcommand{\vp}{\varphi}

\newcommand{\mT}{\mathring{T}}

\newcommand{\zf}{\zeta\left(\frac{1}{2}+it\right)}

\newtheorem{theorem}{Theorem}

\theoremstyle{definition}

\theoremstyle{remark}
\newtheorem{remark}[]{Remark}

\numberwithin{equation}{section}

\begin{document}

\title{Jacob's ladders, Gram's sequence and some nonlinear integral equations connected with the functions $J_\nu(x)$ and $|\zf|^4$}

\author{Jan Moser}

\address{Department of Mathematical Analysis and Numerical Mathematics, Comenius University, Mlynska Dolina M105, 842 48 Bratislava, SLOVAKIA}

\email{jan.mozer@fmph.uniba.sk}


\begin{abstract}
It is shown in this paper that the Jacob's ladder is the asymptotic solution to the new nonlinear integral equations which correspond to the
functions $J_\nu(x)$ and $|\zf|^4$.
\end{abstract}

\maketitle

\section{The first result: the nonlinear integral equation connected with the Bessel's functions}

\subsection{}

In this paper we obtain some new properties of the signal
\bdis
Z(t)=e^{i\vartheta(t)}\zf
\edis
that is generated by the Riemann zeta-function, where
\bdis
\vartheta(t)=-\frac t2\ln\pi+\text{Im}\ln\Gamma\left(\frac 14+i\frac t2\right)=\frac t2\ln\frac{t}{2\pi}-\frac t2-\frac{\pi}{8}+
\mcal{O}\left(\frac{1}{t}\right) .
\edis
Let us remind that
\bdis
\tilde{Z}^2(t)=\frac{{\rm d}\vp_1(t)}{{\rm d}t},\ \vp_1(t)=\frac 12\vp(t)
\edis
where
\be \label{1.1}
\tilde{Z}^2(t)=\frac{Z^2(t)}{2\Phi^\prime_\vp[\vp(t)]}=\frac{\left|\zf\right|^2}{\left\{ 1+\mcal{O}\left(\frac{\ln\ln t}{\ln t}\right)\right\}\ln t}
\ee
(see \cite{1}, (3.9); \cite{2}, (1.3); \cite{7}, (1.1), (3.1), (3.2)), and $\vp(t)$ is the Jacob's ladder, i.e. a solution of the nonlinear integral
equation (see \cite{1})
\bdis
\int_0^{\mu[x(T)]}Z^2(t)e^{-\frac{2}{x(T)}t}{\rm d}t=\int_0^TZ^2(t){\rm d}t .
\edis

\subsection{}

The Gram's sequence $\{ t_\nu\}$ is defined by the equation
\bdis
\vartheta(t_\nu)=\pi\nu,\ \nu=1,2,\dots
\edis
where (see \cite{20}, p. 102)
\be \label{1.2}
t_{\nu+1}-t_\nu=\frac{2\pi}{\ln t_\nu}+\frac{2\pi\ln 2\pi}{\ln^2 t_\nu}+\mcal{O}\left(\frac{1}{\ln^3t_\nu}\right) .
\ee
The following theorem holds true.

\begin{theorem}
Every Jacob's ladder $\vp_1(t)=\frac{1}{2}\vp(t)$ where $\vp(t)$ is the exact solution to the nonlinear integral equation
\bdis
\int_0^{\mu[x(T)]}Z^2(t)e^{-\frac{2}{x(T)}t}{\rm d}t=\int_0^TZ^2(t){\rm d}t
\edis
is the asymptotic solution of the following nonlinear integral equation
\be \label{1.3}
\int_{x^{-1}(t_\nu)}^{x^{-1}(t_{\nu+1})}J_1[x(t)]\left|\zf\right|^2{\rm d}t=\frac{2\sqrt{2\pi}}{\sqrt{t_\nu}}\sin\left(t_\nu-\frac{\pi}{4}\right)
\ee
for every sufficiently big $t_\nu$ that fulfils the conditions
\be \label{1.4}
\begin{split}
& [t_\nu,t_{\nu+1}]\subset [\mu_n^{(1)},\mu_{n+1}^{(1)}], \\
& [t_\nu,t_{\nu+1}]\cap [k\pi-\epsilon,k\pi+\epsilon]=\emptyset,\ \nu,k\in\mbb{N},\ \nu\to\infty
\end{split}
\ee
where $J_1(\mu_n^{(1)})=0,\ n=1,2,\dots $, i.e. the following asymptotic formula
\be \label{1.5}
\int_{\vp_1^{-1}(t_\nu)}^{\vp_1^{-1}(t_{\nu+1})}J_1[\vp_1(t)]\left|\zf\right|^2{\rm d}t\sim\frac{2\sqrt{2\pi}}{\sqrt{t_\nu}}
\sin\left(t_\nu-\frac{\pi}{4}\right)
\ee
holds true.
\end{theorem}

\begin{remark}
Since
\be \label{1.6}
\mu_{n+1}^{(1)}-\mu_n^{(1)}\sim\pi,\ n\to\infty
\ee
then the number $N_{\nu,n}$ of the intervals $[t_\nu,t_{\nu+1}]$, for which the first condition in (\ref{1.4}) is fulfilled, is given by 
the asymptotic formula 
\bdis
N_{\nu,n}\sim\frac 12\ln t_\nu,\ t_\nu\to\infty , 
\edis
((\ref{1.2}), (\ref{1.6})). 
\end{remark}

This paper is a continuation of the series \cite{1} - \cite{19}.

\section{The second result: some nonlinear integral equation connected with the function $|\zf|^4$}

Let us remind that the Jacob's ladder $\vp_2(T)$ of the second order is a solution to the nonlinear integral equation

\be \label{2.1}
\int_0^{\mu[x(T)]} Z^4(t)e^{-\frac{t}{x(T)}}{\rm d}t=\int_0^T Z^4(t){\rm d}t
\ee
(see \cite{8}). In this case the following asymptotic formula (see \cite{8}, (1.5))
\be \label{2.2}
\begin{split}
& \int_{\vp_1^{-1}(T)}^{\vp_1^{-1}(T+U)}\left|\zeta\left(\frac 12+i\vp_2(t)\right)\right|^4\left|\zf\right|^4{\rm d}t\sim \\
& \sim\frac{1}{4\pi^4}U\ln^8T,\ U=T^{13/14+2\epsilon},\ T\to\infty
\end{split}
\ee
holds true.

\begin{remark}
The small improvements of the exponent $\frac{13}{14}$ that are of the type $\frac{13}{14}\rightarrow \frac{8}{9}\rightarrow \dots$ are
irrelevant in this question.
\end{remark}

Next, similarly to the Theorem 1, the following theorem holds true.

\begin{theorem}
Every Jacob's ladder of the second order $\vp_2(t)$, i.e. the (exact) solution to the nonlinear integral equation (\ref{2.1}) is the asymptotic
solution of the nonlinear integral equation
\be \label{2.3}
\int_{x^{-1}(T)}^{x^{-1}(T+U)}\left|\zeta\left(\frac 12+ix(t)\right)\right|^4\left|\zf\right|^4{\rm d}t=\frac{1}{4\pi^4}U\ln^8T,\ T\to\infty ,
\ee
(comp. (\ref{2.2})).
\end{theorem}

\begin{remark}
There are the fixed-point methods and other methods of the functional analysis used to study the nonlinear equations. What can be obtained by using
these methods in the case of the nonlinear integral equations (\ref{1.3}), (\ref{2.3})?
\end{remark}

\section{Proof of the Theorem 1}

\subsection{}

Let us remind that the following lemma holds true (see \cite{6}, (2.5); \cite{7}, (3.3)): for every integrable function (in the Lebesgue sense)
$f(x),\ x\in [\vp_1(T),\vp_1(T+U)]$ we have
\be \label{3.1}
\int_T^{T+U}f[\vp_1(t)]\tilde{Z}^2(t){\rm d}t=\int_{\vp_1(T)}^{\vp_1(T+U)}f(x){\rm d}x,\ U\in \left(\left. 0,\frac{T}{\ln T}\right]\right.
\ee
where
\bdis
t-\vp_1(t)\sim (1-c)\pi(t) ,
\edis
$c$ is the Euler's constant and $\pi(t)$ is the prime-counting function. In the case $\mT=\vp_1^{-1}(T),\ \widering{T+U}=\vp_1^{-1}(T+U)$ we
obtain from (\ref{2.1})
\be \label{3.2}
\int_{\vp_1^{-1}(T)}^{\vp_1^{-1}(T+U)}f[\vp_1(t)]\tilde{Z}^2(t){\rm d}t=\int_T^{T+U}f(x){\rm d}x .
\ee

\subsection{}

By the simple formula
\bdis
\int_0^a J_1(x){\rm d}x=1-J_0(a) ,
\edis
known from the theory of the Bessel's functions, we obtain
\be \label{3.3}
\int_{t_\nu}^{t_{\nu+1}}J_1(x){\rm d}x=J_0(t_\nu)-J_0(t_{\nu+1}) .
\ee

Hence, from (\ref{3.3}) by (\ref{3.2}) the formula
\be \label{3.4}
\int_{\vp_1^{-1}(t_\nu)}^{\vp_1^{-1}(t_{\nu+1})}J_1[\vp_1(t)]\tilde{Z}^2(t){\rm d}t=J_0(t_\nu)-J_0(t_{\nu+1})
\ee
is obtained.

\subsection{}

It is also well-know that

\be \label{3.5}
J_\nu(x)=\sqrt{\frac{2}{\pi x}}\cos\left( x-\nu\frac{\pi}{2}-\frac{\pi}{4}\right)+\mcal{O}\left(\frac{1}{x^{3/2}}\right),\ x\to\infty 
\ee
(the asymptotic formula for $J_\nu(x)$). Since, by the Titchmarsh' formula (\ref{1.2})
\be \label{3.6}
\frac{1}{\sqrt{t_{\nu+1}}}=\frac{1}{\sqrt{t_\nu}}+\mcal{O}\left(\frac{1}{t_\nu^{3/2}\ln t_\nu}\right),
\ee
it follows (see (\ref{3.5}), (\ref{3.6})) that
\be \label{3.7}
\begin{split}
& J_0(t_\nu)-J_0(t_{\nu+1})= \\
& =\sqrt{\frac{2}{\pi t_\nu}}\left\{\cos\left( t_\nu-\frac{\pi}{4}\right)-\cos\left( t_{\nu+1}-\frac{\pi}{4}\right)\right\}+
\mcal{O}\left(\frac{1}{t_\nu^{3/2}}\right) .
\end{split}
\ee
Next, (see (\ref{1.2}))

\be \label{3.8}
\begin{split}
& \cos\left(t_\nu-\frac{\pi}{4}\right)-\cos\left( t_{\nu+1}-\frac{\pi}{4}\right)=2\sin\frac{t_{\nu+1}-t_\nu}{2}
\sin\left(\frac{t_{\nu+1}+t_\nu}{2}-\frac{\pi}{4}\right)= \\
& =2\sin\frac{t_{\nu+1}-t_\nu}{2}\sin\left(\frac{t_{\nu+1}-t_\nu}{2}+t_\nu-\frac{\pi}{4}\right) = \\
& =2\sin^2\frac{t_{\nu+1}-t_\nu}{2}\cos\left(t_\nu-\frac{\pi}{4}\right)+\sin(t_{\nu+1}-t_\nu)\sin\left( t_\nu-\frac{\pi}{4}\right) = \\
& =\frac{2\pi}{\ln t_\nu}\sin\left( t_\nu-\frac{\pi}{4}\right)+\mcal{O}\left(\frac{1}{\ln^2 t_\nu}\right) .
\end{split}
\ee
Hence, from (\ref{3.4}) by (\ref{3.7}), (\ref{3.8}) the asymptotic formula
\be \label{3.9}
\begin{split}
& \int_{\vp_1^{-1}(t_\nu)}^{\vp_1^{-1}(t_{\nu+1})}J_1[\vp_1(t)]\tilde{Z}^2(t){\rm d}t= \\
& =\frac{2\sqrt{2\pi}}{\sqrt{t_\nu}\ln t_\nu}\sin\left( t_\nu-\frac{\pi}{4}\right)+\mcal{O}\left(\frac{1}{\sqrt{t_\nu}\ln^2 t_\nu}\right)
\end{split}
\ee
follows if the second condition in (\ref{1.4}) is fulfilled. Then from (\ref{3.9}) by the mean-value theorem (see (\ref{1.1}), (\ref{1.4}) and
\cite{19}, (3.3)) we obtain

\bdis
\begin{split}
& \int_{\vp_1^{-1}(t_\nu)}^{\vp_1^{-1}(t_{\nu+1})}J_1[\vp_1(t)]\left|\zf\right|^2{\rm d}t= \\
& =\frac{2\sqrt{2\pi}}{\sqrt{t_\nu}}\sin\left( t_\nu-\frac{\pi}{4}\right)+\mcal{O}\left(\frac{\ln\ln t_\nu}{\sqrt{t_\nu}\ln t_\nu}\right) ,
\end{split}
\edis
i.e. the formula (\ref{1.5}) holds true.

\thanks{I would like to thank Michal Demetrian for helping me with the electronic version of this work.}


\begin{thebibliography}{29}
\bibitem{16}
J. Moser, `Jacob's ladders and the nonlocal interaction of the function $Z^2(t)$ with the function $\tilde{Z}^2(t)$ on the distance $\sim (1-c)\pi(t)$
for the collections of disconnected sets`, (2010), arXiv: 1007.5147.
\end{thebibliography}
\end{document}